\documentclass{article}
\usepackage{amssymb}	% Allows use of AMS's mathematical symbols
\usepackage{amsmath}    % Allows use of AMS's mathematical commands
\usepackage{longtable}  % Allows use of tables longer than one page
\usepackage{float}      % place the figures in right place
\usepackage{latexsym}	% Allows use of old latex symbols

\usepackage{epsfig}
\newtheorem{algorithm}{Algorithm}[section]

\newtheorem{theorem}{Theorem}[section]

\newtheorem{remark}{Remark}[section]

\newcommand{\qed}{\nobreak \ifvmode \relax \else \ifdim\lastskip<1.5em \hskip-\lastskip \hskip1.5em plus0em minus0.5em \fi \nobreak \vrule height0.75em width0.5em depth0.25em\fi}

\def\R{{\bf R}}

\def\T{{\rm T}}

\title{A Globally and Quadratically Convergent Algorithm with Efficient Implementation
for Unconstrained Optimization\thanks{Yaguang Yang is with the Office of Research, 
US NRC, 21 Church Street, Rockville, Maryland 20850, The United States. Email: 
yaguang.yang@verizon.net.}}
\author {Yaguang Yang}
\date{\today}

\newpage
\begin{document}

\maketitle    % This command generates the title.

\begin{abstract}
In this paper, an efficient modified Newton type algorithm is proposed for nonlinear
unconstrianed optimization problems. The modified Hessian is a convex combination 
of the identity matrix (for steepest descent algorithm) and the Hessian matrix 
(for Newton algorithm). The coefficients of the convex combination
are dynamically chosen in every iteration. The algorithm is proved to be
globally and quadratically convergent for (convex and nonconvex) nonlinear
functions. Efficient implementation
is described. Numerical test on widely used CUTE test problems
is conducted for the new algorithm. The test results are compared with
those obtained by MATLAB optimization toolbox function {\tt fminunc}.
The test results are also compared with those obtained by some established 
and state-of-the-art algorithms, such as a
limited memory BFGS, a descent and conjugate gradient algorithm, and 
a limited memory and descent conjugate gradient algorithm.
The comparisons show that the new algorithm is promising.
\end{abstract}


\newpage
\section{ Introduction}
Newton type algorithm is attractive due to its fast convergence 
rate \cite{bertse96}. In non-convex case, Newton algorithm may
not be globally convergent, therefore, various modified Newton algorithms have
been proposed, for example \cite{fgm95} \cite{gm74}. The idea is 
to add a positive diagonal matrix to the Hessian matrix so that the
modified Hessian is positive definite and the modified algorithms become globally 
convergent, which is similar to the idea of Levenberg-Marquardt method 
studied in \cite{ysa89}. However, for the iterates far
away from the solution set, the added diagonal matrix may be very
large. This may lead to the poor condition number of the modified Hessian,
generate a very small step, and prevent the iterates from quickly
moving to the solution set \cite{fy05}. 

In this paper, we propose a slightly different modified Newton algorithm. The modified Hessian
is a convex combination of the Hessian (for Newton algorithm) and the identity 
matrix (for steepest descent algorithm). Therefore, the condition number of 
the modified Hessian is well controlled, and the steepest descent algorithm 
and Newton algorithm are special cases of the proposed algorithm. We will show
that the proposed algorithm has merits of both the steepest descent algorithm and 
the Newton algorithm, i.e., the algorithm is globally and quadratically convergent. 
We will also show that the algorithm can be implemented in an efficient way,
using the optimization techniques on Riemannian manifolds proposed in
\cite{eas98}, \cite{smith93}, \cite{yyang96}, and \cite{yang06a}. 
Numerical test for the new algorithm is conducted for the widely used nonlinear 
optimization test problem set CUTE downloaded from \cite{web}. The test
results are compared with those obtained by MATLAB optimization toolbox 
function {\tt fminunc}. The test results are also compared with those obtained by 
some established and state-of-the-art algorithms, such as 
limited memory BFGS \cite{nocedal80}, a descent and conjugate gradient algorithm
\cite{hz05}, and a limited memory and descent conjugate gradient algorithm \cite{hz12}.
The comparison shows that the new algorithm is promising.

The rest paper is organized as follows. Section 2 proposes the modified Newton
algorithm and provides the convergence results. Section 3 discusses an efficient
implementation involving calculations of the maximum and minimum eigenvalues of the modified
Hessian matrix. Section 4 presents numerical test results. The last section 
summarizes the main result of this paper.

\section{Modified Newton Method}
Our objective is to minimize a multi-variable nonlinear (convex or non-convex) function
\begin{equation}
\min f(x),
\end{equation}
where $f$ is twice differentiable. Throughout the paper, we define by $g(x)$ or
simply by $g$ the gradient of $f(x)$, by $H(x)$ or simply by $H$ the Hessian of $f(x)$, 
by $\lambda_{max}H(x)$ or simply $\lambda_{max}(H)$ the maximum eigenvalue of $H(x)$, 
by $\lambda_{min}H(x)$ or simply $\lambda_{min}(H)$ the minimum eigenvalue 
of $H(x)$. Assuming that $\bar{x}$ is a local minimizer, we make the following assumptions 
in our convergence analysis.
\newline
\newline
{\bf Assumptions:}
\begin{itemize}
\item[1.] $g(\bar{x})=0$.
\item[2.] The gradient $g(x)$ is Lipschitz continuous, i.e., there exists a constant
$L>0$ such that for any $x$ and $y$, 
\begin{equation}
\|g(x)-g(y)\| \le L \|x-y\|.
\end{equation}
\item[3.] There are small positive numbers $\delta >0$, $\eta>0$, and a large positive
number $\Delta \ge 1$, and a neighborhood of $\bar{x}$, defined by 
${\cal N}(\bar{x})=\{x: \|g({x})-g(\bar{x}) \| \le \eta \}$, such that for all
$x \in {\cal N}(\bar{x})$, $\lambda_{min}(H({x})) \ge \delta >0$
and $\lambda_{max}(H)/\lambda_{min}(H) \le \Delta$.
\end{itemize}
Assumptions 1 is standard, i.e., $\bar{x}$ meets the first order necessary condition.
If the gradient is Lipschitz continuous as defined in Assumption 2, then ${\cal N}(\bar{x})$ 
is well defined. Assumption 3 indicates that for all $x \in {\cal N}(\bar{x})$, a 
strong second order sufficient condition 
holds, and the condition number of Hessian is bounded which is equivalent to 
$\lambda_{max}(H)<\infty$ given $\lambda_{min}(H({x})) \ge \delta$.

In the remaining discussion, we will use subscript $k$ for the $k$th iteration.
The idea of the proposed algorithm is to search optimizers 
along a direction $d_k$ that satisfies 
\begin{equation}
(\gamma_k I +(1-\gamma_k) H(x_k)) d_k = B_kd_k = - g(x_k),
\label{direction}
\end{equation}
where $\gamma_k \in [0,1]$ will be carefully selected in every iteration.
Clearly, the modified Hessian is a convex combination of the identity matrix
for steepest descent algorithm and the Hessian for the Newton algorithm. 
When $\gamma_k=1$, the algorithm reduces to the steepest descent
algorithm; when $\gamma_k=0$, the algorithm reduces to the Newton algorithm. 
We will focus on the selection of $\gamma_k$, and we will prove the global 
and quadratic convergence of the proposed algorithm. The convergence properties 
are directly related to the goodness of the search direction and step length,
%\cite{nocedal92}
which in turn decide the selection criteria of $\gamma_k$. The quality of the 
search direction is measured by
\begin{equation}
\cos(\theta_k)=-\frac{g_k^Td_k}{\|g_k\| \|d_k\|},
\label{cos}
\end{equation}
which should be bounded below from zero in all iterations. 
A good step length $\alpha_k$ should satisfy the following Wolfe condition.
\begin{subequations}
\begin{align}
f(x_k+\alpha_k d_k) \le f(x_k)+\sigma_1 \alpha_kg_k^Td_k, \\
g(x_k+\alpha_k d_k) \ge \sigma_2 g_k^Td_k,
\end{align}
\end{subequations}
where $0 < \sigma_1 < \sigma_2<1$. The existence of Wolfe condition is 
established in \cite{wolfe69, wolfe71}.
The proposed algorithm is given as follows.
\begin{algorithm} {\bf Modified Newton} 
\\*
Data: $0<\delta$, and $1 \le \Delta<\infty$, initial $x_0$.
\newline
{\bf for} k=0,1,2,...
\begin{itemize}
\item[] Calculate gradient $g(x_k)$.
\item[] Calculate Hessian $H(x_k)$, select $\gamma_k$, and calculate $d_k$ from (\ref{direction}).
\item[] Select $\alpha_k$ and set $x_{k+1}=x_k+\alpha_k d_k$.
\end{itemize}
{\bf end}
\label{newAlgo}
\end{algorithm}

\begin{remark}
An algorithm that finds, in finite steps, a point satisfying Wolfe condition is
given in \cite{more90}. Therefore, the selection of $\alpha_k$ will not be discussed 
in this paper. 
\end{remark}

We will use an important global convergence result given by Zoutendijk
\cite{zoutendijk70} which can be stated as follows.

\begin{theorem}
Suppose that $f$ is bounded below in $\R^n$ and that $f$ is continuously 
twice differentiable
in a neighborhood ${\cal M}$ of the level set ${\cal L}=\{ x: f(x) \le f(x_0) \}$.
Assume that the gradient is Lipschitz continuous 
for all $x, y \in {\cal M}$. Assume further that $d_k$ is a descent direction
and $\alpha_k$ satisfies the Wolfe condition. 
Then
\begin{equation}
\sum_{k \ge 0} \cos^2(\theta_k) \| g_k \|^2 < \infty.
\end{equation}
\hfill \qed
\label{Zoutendijk}
\end{theorem}

Zoutendijk theorem indicates that if for all $k \ge 0$, $d_k$ is a descent direction;
and for a constant $C$, $\cos(\theta_k) \ge C >0$, then the algorithm
is globally convergent because $\lim_{k \rightarrow \infty} \|g_k\|=0$. 
To assure that $d_k$ is a descent direction, $B_k$ should be strictly positive.
This can be achieved by setting
\begin{equation}
\gamma_{k} +(1-\gamma_{k}) \lambda_{min}(H_k) \ge \delta,
\label{condition}
\end{equation}
which is equivalent to
\begin{equation}
\gamma_{k} (1-\lambda_{min}(H_k)) \ge \delta-\lambda_{min}(H_k).
\label{gamma1}
\end{equation}
Therefore, we set
\begin{equation}
\gamma_{k} = \left\{
\begin{array}{ll}
\frac{\delta-\lambda_{min}(H_k)}{1-\lambda_{min}(H_k)} & \mbox{if $\lambda_{min}(H_k)< \delta$} \\
0 & \mbox{if $\lambda_{min}(H_k) \ge \delta$}.
\end{array}
\right.
\label{cond1}
\end{equation}
In view of (\ref{direction}) and (\ref{cos}), it is clear that if 
\begin{equation}
\| B_k \| \|B_k^{-1} \| \le \Delta,
\label{gCond}
\end{equation}
where $1 \le \Delta < \infty$, then $\cos(\theta_k) \ge 1/\Delta =C >0$. Therefore, in view of 
Theorem \ref{Zoutendijk}, to achieve the global
convergence, from (\ref{direction}) and (\ref{gCond}), $\Delta$ should meet the following condition
\begin{equation}
\frac{\gamma_{k} +(1-\gamma_{k}) \lambda_{max}(H_k)}{\gamma_{k} +(1-\gamma_{k}) \lambda_{min}(H_k)} \le \Delta.
\end{equation}
Using (\ref{condition}) and $\gamma_{k} +(1-\gamma_{k}) \lambda_{min}(H_k)>0$, we have
\begin{equation}
({\Delta -1+\lambda_{max}(H_k)-\lambda_{min}(H_k) \Delta})\gamma_{k} \ge 
{\lambda_{max}(H_k)-\lambda_{min}(H_k) \Delta}.
\label{gamma2}
\end{equation}
Since ${\Delta -1+\lambda_{max}(H_k)-\lambda_{min}(H_k) \Delta} \ge
{\lambda_{max}(H_k)-\lambda_{min}(H_k) \Delta}$, we should select
\begin{equation}
\gamma_{k} = \left\{
\begin{array}{ll}
0 & \mbox{if ${\lambda_{max}(H_k)} \le {\lambda_{min}(H_k)} \Delta$}  \\
\frac{\lambda_{max}(H_k)-\lambda_{min}(H_k) \Delta}
{\Delta -1+\lambda_{max}(H_k)-\lambda_{min}(H_k) \Delta} 
& \mbox{if ${\lambda_{max}(H_k)} > {\lambda_{min}(H_k)} \Delta$}.
\end{array}
\right.
\label{cond2}
\end{equation}
Combining (\ref{cond1}) and (\ref{cond2}) yields
\begin{equation}
\gamma_k = \left\{
\begin{array}{ll}
0 & \mbox{if $\lambda_{min}(H_k) \ge \delta$ and ${\lambda_{max}(H_k)} \le {\lambda_{min}(H_k)}\Delta$} \\
a_k=\frac{\delta-\lambda{min}(H_k)}{1-\lambda{min}(H_k)} 
   & \mbox{if $\lambda_{min}(H_k) < \delta$ and 
   ${\lambda_{max}(H_k)} \le {\lambda_{min}(H_k)}\Delta$}  \\
b_k=\frac{\lambda_{max}(H_k)-\lambda_{min}(H_k) \Delta}
{\Delta -1+\lambda_{max}(H_k)-\lambda_{min}(H_k) \Delta} 
    &  \mbox{if $\lambda_{min}(H_k) \ge \delta$ and
    	${\lambda_{max}(H_k)} > {\lambda_{min}(H_k)}\Delta $}   \\
\max \{ a_k,b_k \} & \mbox{else}
\end{array}
\right.
\label{gamma}
\end{equation}
It is clear to see from the selection of $\gamma_k$ that (\ref{gamma1}) and (\ref{gamma2}) hold.
This means that the conditions of Theorem \ref{Zoutendijk} hold. Therefore, Algorithm \ref{newAlgo} 
is globally convergent.
 
Since Algorithm \ref{newAlgo} is globally convergent in the sense that
$\lim_{k \rightarrow \infty} \|g_k\|=0$, there exists an $\eta>0$ such that for
$k$ large enough, $\| g(x_k) \| \le \eta$; from Assumption 3, $\lambda_{min}(H({x}_k)) \ge \delta >0$
and $\lambda_{max}(H_k)/\lambda_{min}(H_k) \le \Delta$. From (\ref{gamma}), 
$\gamma_k=0$ for all $k$ large enough, i.e., Algorithm \ref{newAlgo} reduces to
Newton algorithm. Therefore, the proposed algorithm is quadratic convergent. We summarize the
main result of this paper as the following 

\begin{theorem}
Suppose that $f$ is bounded below in $\R^n$ and that $f$ is continuously 
twice differentiable
in a neighborhood ${\cal M}$ of the level set ${\cal L}=\{ x: f(x) \le f(x_0) \}$.
Assume that the gradient is Lipschitz continuous
for all $x, y \in {\cal M}$. Assume further that $d_k$ is defined as in (\ref{direction})
with $\gamma_k$ being selected as in (\ref{gamma}) and $\alpha_k$ satisfies the Wolfe condition. 
Then Algorithm \ref{newAlgo} is globally convergent. Moreover, if the
convergent point $\bar{x}$ satisfies Assumption 3, then Algorithm \ref{newAlgo} converges
to $\bar{x}$ in quadratic rate.
\hfill \qed
\end{theorem}

%\begin{remark}
%When $x_k$ is close to a local maximizer where $\lambda_{max}(H_k) < 0$. The Newton algorithm
%may converge to the maximizer because it directly solve $g(x)=0$ which is satisfied by the local
%maximizer. But the modified Newton algorithm will search along the steepest descent direction 
%for minimizers because $\gamma_k=1$ is selected.
%\end{remark}

\begin{remark}
Since $B_k$ is positive definite, Cholesky factorization exist, (\ref{direction}) can be solved 
efficiently. Furthermore, if $H_k$ is sparse, (\ref{direction}) can be solved using techniques 
for a sparse matrix.
\end{remark}

\section{Implementation Consideration}
To implement Algorithm \ref{newAlgo} for practical use, we need to consider several issues.

\subsection{Termination}
First, we need to have a termination rule in Algorithm \ref{newAlgo}. This rule is checked at the end of 
Step 1. For $0<\epsilon$, if $\| g(x_k) \| < \epsilon$ or $\| g(x_k) \|_{\infty} < \epsilon$, stop.

\subsection{Computation of extreme eigenvalues}
The most significant computation in the proposed algorithm is the selection of $\gamma_k$, which 
involves the computation of $\lambda_{max}(H_k)$ and $\lambda_{min}(H_k)$ for the symmetric
matrix $H$. There are general algorithms to compute eigenvalues and eigenvectors for a
symmetric matrix \cite{gv89}. However, there are much more efficient algorithms for extreme 
eigenvalues for a symmetric matrix, which is equivalent to find the solution of Rayleigh Quotient
\begin{equation}
\lambda_{max}(H_k) = \max_{\| x \|=1} x^{\T}H_k x =\max \frac{x^{\T}H_k x}{x^{\T} x},
\label{maxH}
\end{equation}
\begin{equation}
\lambda_{min}(H_k) = \min_{\| x \|=1} x^{\T}H_k x =\min \frac{x^{\T}H_k x}{x^{\T} x}.
\label{minH}
\end{equation}
It is well-known that there are cubically convergent algorithms to find the solution of 
Rayleigh Quotient \cite{eas98}. In our opinion, the most efficient methods are the 
conjugate gradient optimization algorithm on Riemannian manifold proposed by Smith \cite{smith93},
and the Armijo-Newton optimization algorithm on Riemannian manifold proposed by Yang \cite{yyang96}
\cite{yang06a}. Both methods make fully use of the geometry of the unit sphere ($\| x \|=1$) 
and search the solution along the arc defined by geodesics over the unit sphere. Armijo-Newton 
algorithm may converge faster, but it may converge to an internal eigenvalue rather than an 
extreme eigenvalue. Conjugate gradient optimization algorithm may also converge to an internal 
eigenvalue, but the chance is much smaller and a small perturbation may lead the iterate
to converge to the desired extreme eigenvalues. Let $x$ be on unit sphere and $\rho(x)=x^{\T}Hx$. 
For vector $v$ in tangent space at $x$, $\tau v$ denote parallelism of $v$ along the geodesic 
defined by a unit length vector $q$ in tangent space at $x$, it is shown in \cite{smith93}
\begin{equation}
\tau v = v-(v^{\T}q)(x\sin(t)+q(1-\cos(t)).
\label{parallel}
\end{equation}
To find the maximum eigenvalue of $H$ defined in (\ref{maxH}), the conjugate 
gradient algorithm proposed in \cite{smith93} is stated as follows 
(with very minor but important modification presented in bold font).

\begin{algorithm} {\bf Conjugate gradient (CG) for maximum eigenvalue} 
\\*
Data: $0<\epsilon$, initial $x_0$ with $\| x_0 \|=1$, 
$G_0=Q_0=(H-\rho(x_0) I)x_0$.
\newline
{\bf for} k=0,1,2,...
\begin{itemize}
\item[] Calculate $c$, $s$, and $v=1-c$, such that $\rho(x_k c+q_ks)$ is maximized,
where $c^2+s^2=1$, $q_k=\frac{Q_k}{\| Q_k \|}$. This can be accomplished by 
geodesic maximization and the formula is given by
\begin{equation}
\left\{
\begin{array}{lll}
c=\sqrt{\frac{1}{2}\left(1+\frac{b}{r}\right)}, & s=\frac{a}{2rc}, & \mbox{if $b \ge 0$} \\
s=\sqrt{\frac{1}{2}\left(1-\frac{b}{r}\right)}, & c=\frac{a}{2rs}, & \mbox{if $b \le 0$}
\end{array}
\right.
\label{csMax}
\end{equation}
where $a=2x_k^{\T}Hq_k$, $b=x_k^{\T}Hx_k-q_k^{\T}Hq_k$, and $r=\sqrt{a^2+b^2}$.
\item[] Set $x_{k+1}=x_k c+q_k s$, ${\bf x_{k+1}=\frac{x_{k+1}}{\| x_{k+1} \|}}$, 
$\tau Q_k=Q_k c-x_k \|Q_k\|s$, and
$\tau G_k=G_k-(q_k^{\T}G_k)(x_k s+q_k v)$.
\item[] Set $G_{k+1}=(H-\rho(x_{k+1})I)x_{k+1}$, $Q_{k+1}=G_{k+1}+\mu_k \tau Q_k$, 
where $\mu_k=\frac{(G_{k+1}-\tau G_k)^{\T}G_{k+1}}{G_k^{\T}Q_k}$. 
\item[] Set ${\bf Q_{k+1}=(I-x_{k+1}x_{k+1}^{\T})Q_{k+1} }$.
\item[] If $k \equiv n-1$, set $Q_{k+1}={\bf (I-x_{k+1}x_{k+1}^{\T})}G_{k+1}$.
\end{itemize}
{\bf end}
\label{CG}
\end{algorithm}

\begin{remark}
$Q_{k+1}$ should be on tangent space at $x_{k+1}$. But numerical error may change $Q_{k+1}$
slightly. Therefore, the projection is necessary to bring $Q_{k+1}$ back to the tangent 
space at $x_{k+1}$. Similar changes are made to ensure the unit length of $x_k$.
With these minor changes, the CG algorithm is much more stable and the 
observed convergence rate is faster than the one reported in \cite{smith93}.
\end{remark}

\begin{remark}
To search for the minimum eigenvalue of (\ref{minH}), (\ref{csMax}) is replaced by
\begin{equation}
\left\{
\begin{array}{lll}
c=\sqrt{\frac{1}{2}\left(1-\frac{b}{r}\right)}, & s=\frac{a}{2rc}, & \mbox{if $b \ge 0$} \\
s=\sqrt{\frac{1}{2}\left(1+\frac{b}{r}\right)}, & c=\frac{a}{2rs}, & \mbox{if $b \le 0$}
\end{array}
\right.
\label{csMin}
\end{equation}
which is obtained by minimizing $\rho(x_k c+q_ks)$ under the constraint $c^2+s^2=1$.
\end{remark}

\begin{remark}
Each iteration of Algorithm \ref{CG} involves only matrix
and vector multiplications, the cost ${\cal O}(n^2)$ is very low. Our experience 
shows that it needs only a few iterations to converge to the extreme eigenvalues.
\end{remark}

\begin{remark}
If $H$ is sparse, Algorithm \ref{CG} will be very efficient.
\end{remark}

\subsection{The implemented modified Newton algorithm}

The implemented modified Newton algorithm is as follows.

\begin{algorithm} {\bf Modified Newton} 
\\*
Data: $0<\delta$, $0<\epsilon$, and $1 \le \Delta<\infty$, initial $x_0$ with $\| x_0 \|=1$.
\newline
{\bf for} k=1,2,...
\begin{itemize}
\item[] Calculate gradient $g(x_k)$.
\item[] If $\| g(x_k) \|<\epsilon$ or $\| g(x_k) \|_{\infty}<\epsilon$, stop.
\item[] Calculate Hessian $H(x_k)$. 
\item[] Calculate $\lambda_{max}(H_k)$ and $\lambda_{min}(H_k)$ using Algorithm \ref{CG}.
\item[] Select $\gamma_k$ using (\ref{gamma}), and calculate $d_k$ using (\ref{direction}).
\item[] If $d_k$ is not a descent direction, Algorithm \ref{CG} generates an internal eigevalue. 
	A conventional method will be used to find $\lambda_{max}(H_k)$ and $\lambda_{min}(H_k)$.
	Then, select $\gamma_k$ using (\ref{gamma}), and calculate $d_k$ using (\ref{direction}).
\item[] Select $\alpha_k$ using one dimensional search and set $x_{k+1}=x_k+\alpha_k d_k$.
\end{itemize}
{\bf end}
\label{implementedAlgo}
\end{algorithm}

\begin{remark} It is very rare to use a conventional method to calculate $\lambda_{max}(H_k)$ 
and $\lambda_{min}(H_k)$. But this safeguard is needed in case that Algorithm \ref{CG} 
generates an internal eigevalue.
\end{remark}

\section{Numerical Test}
In this section, we present some test results for both Algorithm \ref{CG} and 
Algorithm \ref{implementedAlgo}. 

\subsection{Test of Algorithm \ref{CG}}
The advantages of Algorithm \ref{CG} have been explained 
in \cite{es96}. We conducted numerical test on some problems to confirm the theoretical 
analysis. For the sake of comparison, we use an example in \cite{ysa89} because it provides
detailed information about the test problem and the results obtained by many other 
algorithms. For this problem, 
\[
H=\left[ \begin{array}{cccc}
r_0 & r_1 & \cdots & r_{15} \\
r_1 & r_0 & \cdots & r_{14} \\
\vdots & \vdots & \cdots & \vdots \\
r_{15} & r_{14} & \cdots & r_0 
\end{array} \right],
\]
where 
$r_0=1.00000000$, $r_1=0.91189350$, $r_2 = 0.75982820$,
$r_3 = 0.59792770$, $r_4 = 0.41953610$, $r_5 = 0.27267350$,
$r_6 = 0.13446390$, $r_7 = 0.00821722$, $r_8 =-0.09794101$,
$r_9 =-0.21197350$, $r_{10}=-0.30446960$, $r_{11}=-0.34471370$,
$r_{12}=-0.34736840$, $r_{13}=-0.32881280$, $r_{14}=-0.29269750$,
$r_{15}=-0.24512650$. The minimum eigenvalue is $\lambda_{min}=0.00325850037049$.
Four methods, namely HE, TJ, FR, and CA, which use formulae derived from \cite{csbd86},  
\cite{fr64}, \cite{tj78}, and \cite{ysa89},, are tested and reported in \cite{ysa89}.
These test results are compared with our test obtained by Algorithm \ref{CG} (CG).
The comparison is presented in Table 1. The result is clearly in favor of Algorithm \ref{CG} (CG).

\begin{table}[H]
\begin{center}
\caption{Simulation results of 5 algorithms for the test problem}
\begin{tabular}{|c|c|c|c|c|}
\hline 
 $x_0$ & \multicolumn{2}{c|} {$(-1,1,-1,\cdots)^{\T}$} & \multicolumn{2}{c|} {$(1,0,\cdots,0)^{\T}$} \\
\hline
Algo      & iter  & $\lambda_{min}$   &  iter  & $\lambda_{min}$    \\ \hline
HE        & 24    & 0.0032585         &  77    & 0.0032586   \\ \hline
TJ        & 26    & 0.0032585  	      &  65    & 0.0032586   \\ \hline
FR        & 17    & 0.0032585  	      &  87    & 0.0032586   \\ \hline
CA        & 32    & 0.0032585  	      &  124   & 0.0032586   \\ \hline
CG        & 10    & 0.0032585  	      &  14    & 0.0032585   \\ \hline
\end{tabular}
\end{center}
\end{table}

\subsection{Test of Algorithm \ref{implementedAlgo} on Rosenbrock function}
Algorithm \ref{implementedAlgo} is implemented in Matlab function {\tt mNewton}. The following
parameters are chosen: $\delta=10^{-8}$, $\Delta=10^{12}$, and $\epsilon=10^{-5}$. 
A test for Algorithm \ref{implementedAlgo} is done for Rosenbrock function given by 
\[
f(x) = 100(x_2-x_1^2)^2+(1-x_1)^2,
\]
with initial point $x_0=[-1.9, 2.0]^{\T}$.
Steepest descent is inefficient in this problem. After 1000 iterations, the iterate is still 
a considerable distance from the minimum point $x^*=[1,1]^{\T}$. BFGS algorithm is significantly
better, after 34 iterations, the iterate terminates at $x=[0.9998, 0.9996]^{\T}$ 
(cf. \cite{matlab}). The new algorithm performs even better, after 24 iterations, the iterate 
terminates at $x=[0.9999, 0.9998]^{\T}$. Similar to BFGS algorithm, the new method is able 
to follow the shape of the valley and converges to the minimum as depicted in Figure 1, where
the contour of the Rosenbrock function, the gradient flow from the initial point to the minimum 
point (in blue line), and all iterates (in red "x")
are plotted. 
\begin{figure}[H]
\centerline{\epsfig{file=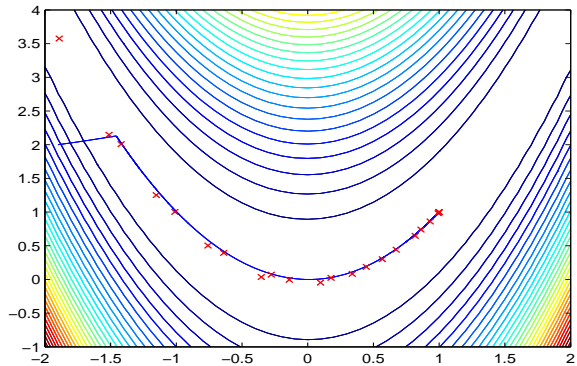,height=5.5cm,width=9cm}}
\caption{New algorithm searches follows the shape of the valley of Rosenbrock function.}
\label{fig:iter1}
\end{figure}

\subsection{Test of Algorithm \ref{implementedAlgo} on CUTE problems}
We also conducted test for both {\tt mNewton} and Matlab optimization toolbox 
function {\tt fminunc} against {CUTE} test problem set.
{\tt fminunc} options are set as
\[
\mbox{options = optimset('MaxFunEvals',1e+20,'MaxIter',5e+5,'TolFun',1e-20, 'TolX',1e-10).}
\]
This setting is selected to ensure that the Matlab function {\tt fminunc} 
will have enough iterations to converge or to fail.
CUTE test problem set is downloaded from Princeton test problem collections \cite{web}. 
Since {CUTE} test set is presented in AMPL mod-files, we first convert AMPL mod-files
into nl-files so that Matlab functions can read the CUTE models, then we use Matlab 
functions {\tt mNewton} and {\tt fminunc} to read the nl-files and solve these test problems.
Because the conversion software which converts mod-files to nl-files is restricted to problems
whose sizes are smaller than $300$,
the test is done for all CUTE unconstrained optimization problems whose sizes 
are less than $300$. The test uses the initial points provided by {CUTE} test problem set, 
we record the calculated objective function values, the norms of the gradients 
at the final points, and the iteration numbers for these testing problems. We present 
the test results in Table 2, and summarize the comparison of the test results as follows:
\begin{itemize}
\item[1.] the modified Newton function {\tt mNewton} converges in all the test problems after
terminate condition $\| g(x_k) \|< 10^{-5}$ is met. But for about $40\%$ of the problems,
Matlab optimization toolbox function {\tt fminunc} does not reduce $\| g(x_k) \|$ to a value
smaller than $0.01$. For these problems, the objective functions obtained by {\tt fminunc} normally 
are not close to the minimum;
\item[2.] for problems that both {\tt mNewton} and {\tt fminunc} converge, {\tt mNewton}
normally uses less iterations than {\tt fminunc} and converges to points with smaller 
$\| g(x_k) \|$ except 2 problems {\tt bard} and {\tt deconvu}.
\end{itemize}

\begin{center}
%\caption{Test results for problems in CUTE \cite{bcgt95}, initial points are given in CUTE}
\begin{longtable}{|c|c|c|c|c|c|c|}
\caption{Test result for problems in CUTE \cite{bcgt95}, initial points are given in CUTE}\\
\hline    
Problem  & iter   &  obj    & gradient  & iter    & obj     & gradient  \\
        & mNewton  &  mNewton  & mNewton     & fminunc & fminunc & fminunc    \\
\hline
%\endfirsthead
%\multicolumn{7}{c}
%{\tablename \thetable\ -- \textit{Continued from previous page}} \\
%\hline
%Problem  & iter   &  obj    & gradient  & iter    & obj     & gradient  \\
%        & mNewton  &  mNewton  & mNewton     & fminunc & fminunc & fminunc    \\
%\hline
%\endhead
%\hline \multicolumn{7}{r}{\textit{Continued on next page}}
%\endfoot
%\hline
%\endlastfoot
%penalty2    & 213 & 6  &   4.46E-09      & & &        \\ \hline
arglina      &  1    & 100.000000       & 0.00000000e-9  & 4  & 100.000000      &  0.00016620       \\ \hline
bard         &  24   & 0.11597205       & 0.96811289e-5  & 20 & 0.00821487      &  0.1158381e-5     \\ \hline
beale        &  6    & 0.00000000e-9    & 0.38908644e-9  & 15 & 0.00000024e-5   &  1.3929429e-5     \\ \hline
%biggs6       &  70   & 0.0.000010e-5    & 0.78652391e-5  & 68 & 0.00047968      &  {\bf 0.01802530} \\ \hline
%box3         &  19   & 0.00000388e-5    & 0.18288127e-5  & 24 & 0.00000388e-5   &  0.2364393e-5     \\ \hline
brkmcc       &  2    & 0.16904267       & 0.61053106e-5  & 5  & 0.16904268      &  0.0454266e-5     \\ \hline
brownal      &  7    & 0.00000000e-7    & 0.26143011e-7  & 16 & 0.00030509e-5   &  0.00010437       \\ \hline
brownbs      &  8    & 0.00000000e-9    & 0.00000000e-9  & 11 & 0.00009308      &  {\bf 15798.5950} \\ \hline
brownden     &  8    & 85822.2016       & 0.00003000e-5  & 32 & 85822.2017      &  {\bf 0.46462733} \\ \hline
chnrosnb     &  46   & 0.00000000e-5    & 0.10455150e-5  & 98 & 30.0583699      &  {\bf 10.1863739} \\ \hline
cliff        &  26   & 0.19978661       & 0.10751025e-6  & 1  & 1.00159994      &  {\bf 1.41477930} \\ \hline
cube         &  28   & 0.00000000e-9    & 0.69055669e-9  & 34 & 0.79877450e-9   &  0.00013409       \\ \hline
deconvu      & 2612  & 0.00242309e-5    & 0.99584075e-5  & 80 & 0.00031582e-3   &  0.1750297e-3     \\ \hline
denschna     &  5    & 0.00000022e-5    & 0.29676520e-5  & 10 & 0.00000005e-5   &  0.1581909e-5     \\ \hline
denschnb     &  5    & 0.00000004e-5    & 0.17646764e-5  & 7  & 0.00000010e-5   &  0.2200204e-5     \\ \hline
denschnc     &  11   & 0.00000000e-9    & 0.17803850e-9  & 21 & 0.00000160e-3   &  0.3262483e-3     \\ \hline
denschnd     &  36   & 0.00126578e-5    & 0.77956675e-5  & 23 & 45.2971677      &  {\bf 84.5851141} \\ \hline
denschnf     &  6    & 0.00000000e-9    & 0.62887898e-9  & 10 & 0.00000002e-3   &  0.1005028e-3     \\ \hline
dixon3dq     &  1    & 0.00000000e-7    & 0.00000000e-7  & 20 & 0.00000014e-5   &  0.3661452e-5     \\ \hline
%djtl         &  14   & -8951.5447       & 0.00002068e-5  & 3  & -8033.8869      &  {\bf 1273.33191} \\ \hline
eigenals     &  22   & 0.00000000e-7    & 0.45589372e-7  & 78 & 0.10928398e-2   &  {\bf 0.10292633} \\ \hline
eigenbls     &  62   & 0.00000000e-6    & 0.32395333e-6  & 91 & 0.34624147      &  {\bf 0.46420894} \\ \hline
engval2      &  13   & 0.00000001e-5    & 0.36978724e-7  & 29 & 0.00003953e-5   &  0.2799583e-3     \\ \hline
%errinros    &  48   & 0.02442067e-5    & 0.51419855e-5  & 92 & 0.00045778      &  {\bf 0.25536229} \\ \hline
%expfit      &  11   & 0.24051059       & 0.23505443e-5  & 12 & 0.24051059      &  0.2263669e-5     \\ \hline
extrosnb     &  1    & 0.00000000       & 0.00000000     & 1  & 0.00000000      &  0.00000000       \\ \hline
fletcbv2     &  1    & -0.5140067       & 0.50699056e-5  & 98 & -0.5140067      &  0.1087190e-4     \\ \hline
fletchcr     &  12   & 0.00000000e-7    & 0.12606909e-7  & 63 & 68.128920       &  {\bf 160.987949} \\ \hline
genhumps     &  52   & 0.00000003e-7    & 0.29148635e-7  & 59 & 0.00044932e-3   &  0.3167733e-3     \\ \hline
%growthls    &  1    & {\bf 3542.14903} & 0.00000000e-5  & 12 & 12.4523620      &  {\bf 0.05809928} \\ \hline
hairy        &  19   & 20.0000000       & 0.00065611e-5  & 22 & 20.000000000    &  0.3810773e-4     \\ \hline
%hatfldd     &  23   & 0.00661522e-5    & 0.18530933e-5  & 19 & 0.00661573e-5   &  0.2355759e-5     \\ \hline
%hatflde     &  28   & 0.04434795e-5    & 0.53212673e-5  & 9  & 0.06210367e-5   &  0.7970000e-5     \\ \hline
heart6ls     & 375   & 0.00000000e-5    & 0.29136580e-5  & 53 & 0.63188192      &  {\bf 71.9382548} \\ \hline
helix        &  13   & 0.00000000e-9    & 0.31818245e-9  & 29 & 0.00000226e-5   &  0.4196860e-4     \\ \hline
hilberta     &  1    & 0.00001538e-7    & 0.92172479e-7  & 35 & 0.02289322e-5   &  0.3263435e-5     \\ \hline
hilbertb     &  1    & 0.0000004e-20    & 0.1267079e-12  & 6  & 0.00000021e-5   &  0.6542441e-5     \\ \hline
himmelbb     &  7    & 0.0000001e-13    & 0.13251887e-6  & 6  & 0.00001462      &  0.0012511        \\ \hline
%himmelbf    &  7    & 0.00000010e-7    & 0.93377677e-7  & 8  & 0.00000001e-5   &  0.1448585e-5     \\ \hline
%himmelbg    &  7    & 0.00000010e-6    & 0.93377677e-6  & 8  & 0.00000001e-5   &  0.1448585e-5     \\ \hline
himmelbh     &  4    & -1.0000000       & 0.00108475e-6  & 7  & -0.9999999      &  0.2607156e-6     \\ \hline
humps        &  26   & 0.00000379e-6    & 0.27563083e-5  & 25 & 5.42481702      &  {\bf 2.36255440} \\ \hline
jensmp       &  9    & 124.362182       & 0.00480283e-5  & 16 & 124.362182      &  0.2897049e-5     \\ \hline
kowosb       &  10   & 0.30750561e-3    & 0.31930835e-5  & 33 & 0.30750560e-3   &  0.0125375e-5     \\ \hline
loghairy     &  23   & 0.18232155       & 0.00103147e-5  & 11 &{\bf 2.5199616136}&  0.0053770        \\ \hline
mancino      &  4    & 0.00000000e-5    & 0.26436736e-5  & 9  & 0.00220471      &  {\bf 1.22432874} \\ \hline
maratosb     &  7    & -1.0000000       & 0.09342000e-9  & 2  & -0.9997167      &  {\bf 0.03570911} \\ \hline
mexhat       &  4    & -0.0401000       & 0.00000000e-5  & 4  & -0.0400999      &  0.1370395e-4     \\ \hline
%meyer3      & & & & & & \\ \hline
%osbornea    & & & & & & \\ \hline
%osborneb    &  51   & 0.04013773       & 0.24749370e-5  & 76 & 0.04013773      &  0.7884569e-5     \\ \hline
palmer1c     &  6    & 0.09759802       & 0.46161602e-5  & 38 & 16139.4418      &  {\bf 655.015973} \\ \hline
palmer2c     &  1    & 0.01442139       & 0.00107794e-5  & 60 & 98.0867115      &  {\bf 33.4524366} \\ \hline
palmer3c     &  1    & 0.01953763       & 0.00434478e-6  & 56 & 54.3139592      &  {\bf 7.85183915} \\ \hline
palmer4c     &  1    & 0.05031069       & 0.01265948e-6  & 56 & 62.2623173      &  {\bf 6.67991745} \\ \hline
palmer5c     &  1    & 2.12808666       & 0.00000001e-5  & 14 & 2.12808668      &  0.00074844       \\ \hline
palmer6c     &  1    & 0.01638742       & 0.00008202e-5  & 43 & 18.0992853      &  {\bf 0.78517164} \\ \hline
palmer7c     &  1    & 0.60198567       & 0.00120838e-5  & 28 & 56.9098797      &  {\bf 4.02685779} \\ \hline
palmer8c     &  1    & 0.15976806       & 0.00013200e-5  & 49 & 22.4365812      &  {\bf 1.31472249} \\ \hline
%pfit1ls     & & & & & & \\ \hline
%pfit2ls     & & & & & & \\ \hline
%pfit3ls     & & & & & & \\ \hline
%pfit4ls     & & & & & & \\ \hline
powellsq     & 0 & 0 & 0 & 0 & 0 & 0 \\ \hline
rosenbr      &  20   & 0.00000002e-5    & 0.10228263e-5  & 36 & 0.00000283e-5   &  2.6095725e-5     \\ \hline
sineval      &  41   & 0.00000000e-8    & 0.24394083e-8  & 47 & 0.22121569      &  {\bf 1.23159435} \\ \hline
sisser       &  13   & 0.00097741e-5    & 0.51113540e-5  & 11 & 0.15409254e-7   &  0.7282671e-5     \\ \hline
tointqor     &  1    & 1175.47222214    & 0.00000000000  & 40 & 1175.4722221    &  0.0090419e-5     \\ \hline
vardim       &  19   & 0.00000000e-8    & 0.00991963e-8  & 1  & 0.22445009e-6   &  {\bf 0.55115494} \\ \hline
watson       &  13   & 0.15239635e-6    & 0.03339433e-6  & 90 & 0.00105098      &  {\bf 0.48756107} \\ \hline
yfitu        &  36   & 0.00000066e-6    & 0.10418764e-6  & 57 & 0.00439883      &  {\bf 11.8427717} \\ \hline
\end{longtable}
\end{center}

\subsection{Comparison of Algorithm \ref{implementedAlgo} to established and 
state-of-the-art algorithms}

Most of the above problems are also used, for example in \cite{hager}, to test some established and 
state-of-the-art algorithms. 
In \cite{hager}, $145$ CUTEr unconstrained problems are tested against limited memory BFGS algorithm 
\cite{nocedal80} (implemented as {\tt L-BFGS}), a descent and conjugate gradient algorithm \cite{hz05}
(implemented as {\tt CG-Descent 5.3}), and a limited memory descent and conjugate gradient algorithm 
\cite{hz12} (implemented as {\tt L-CG-Descent}). The sizes of most of these test problems are smaller 
than or equal to $300$. The size of the largest test problems in \cite{hager} is $10000$. 
Since our AMPL converion software does not work for problems whose sizes are 
larger than $300$, we compare only problems whose sizes are less than or equal to $300$. 
The test results obtained by algorithms descried in \cite{nocedal80, hz05, hz12} are reported in \cite{hager}.
In this test, we changed the stopping criterion for Algorithm \ref{implementedAlgo} to 
$\| g(x) \|_{\infty} \le 10^{-6}$ for consistency. The test results are listed in Table 3.

\begin{center}
\begin{longtable}{|r|r|r|r|r|r|}
\caption{Comparison of mNewtow, L-CG-Descent, L-BFGS, and CG-Descent 5.3 for problems in CUTE \cite{bcgt95}, initial points are given in CUTE}\\
\hline    
Problem      & size   & methods          & iter   &  obj    & gradient      \\
\hline
arglina      &  200   & mNewtow          &   1     &  1.000e+002             &   3.203e-014  \\ 
             &        & L-CG-Descent     &   1     &  {\bf 2.000e+002}       &   3.384e-008  \\
             &        & L-BFGS           &   1     &  {\bf 2.000e+002}       &   3.384e-008  \\
             &        & CG-Descent 5.3   &   1     &  {\bf 2.000e+002}       &   2.390e-007  \\
\hline
bard         &  3     & mNewtow          &   41     &   1.157e-001      &    9.765e-007           \\ 
             &        & L-CG-Descent     &   16     &   8.215e-003      &    3.673e-009           \\
             &        & L-BFGS           &   16     &   8.215e-003      &    3.673e-009           \\
             &        & CG-Descent 5.3   &   21     &   8.215e-003      &    1.912e-007           \\
\hline
beale        &  2     & mNewtow          &   6      &   4.957e-020      &    2.979e-010           \\ 
             &        & L-CG-Descent     &   15     &   2.727e-015      &    4.499e-008           \\
             &        & L-BFGS           &   15     &   2.727e-015      &    4.499e-008           \\
             &        & CG-Descent 5.3   &   18     &   1.497e-007      &    4.297e-007          \\
\hline
brkmcc       &  2     & mNewtow          &   3     &    1.690e-001     &    5.640e-013           \\ 
             &        & L-CG-Descent     &   5     &    1.690e-001     &    6.220e-008           \\
             &        & L-BFGS           &   5     &    1.690e-001     &    6.220e-008           \\
             &        & CG-Descent 5.3   &   4     &    1.690e-001     &    5.272e-008           \\
\hline
%brownal      &  200   & mNewtow          &   7     &    1.496e-016     &    8.518e-009           \\ 
%             &        & L-CG-Descent     &   9     &    2.704e-018     &    4.540e-008           \\
%             &        & L-BFGS           &   4     &    1.473e-009     &    6.663e-007           \\
%             &        & CG-Descent 5.3   &   12    &    6.562e-011     &    1.392e-007           \\
%\hline
brownbs      &  2     & mNewtow          &   8      &   0.000e+000      &   0.000e+000           \\ 
             &        & L-CG-Descent     &   13     &   0.000e+000      &   0.000e+000            \\
             &        & L-BFGS           &   13     &   0.000e+000      &   0.000e+000            \\
             &        & CG-Descent 5.3   &   16     &   1.972e-031      &   8.882e-010            \\
\hline
brownden     &  4     & mNewtow          &   8      &   8.582e+004      &   3.092e-010            \\ 
             &        & L-CG-Descent     &   16     &   8.582e+004      &   1.282e-007            \\
             &        & L-BFGS           &   16     &   8.582e+004      &   1.282e-007            \\
             &        & CG-Descent 5.3   &   38     &   8.582e+004      &   9.083e-007            \\
\hline
chnrosnb     &  50    & mNewtow          &   46      &  1.885e-014       &  7.155e-007             \\ 
             &        & L-CG-Descent     &   287     &  6.818e-014       &  5.414e-007             \\
             &        & L-BFGS           &   216     &  1.582e-013       &  5.565e-007             \\
             &        & CG-Descent 5.3   &   287     &  6.818e-014       &  5.414e-007             \\
\hline
cliff        &  2     & mNewtow          &   26     &   1.998e-001      &   7.602e-008            \\ 
             &        & L-CG-Descent     &   18     &   1.998e-001      &   2.316e-009            \\
             &        & L-BFGS           &   18     &   1.998e-001      &   2.316e-009            \\
             &        & CG-Descent 5.3   &   19     &   1.998e-001      &   6.352e-008            \\
\hline
cube         &  2     & mNewtow          &   28     &   1.238e-017      &   1.985e-007            \\ 
             &        & L-CG-Descent     &   32     &   1.269e-017      &   1.225e-009            \\
             &        & L-BFGS           &   32     &   1.269e-017      &   1.225e-009            \\
             &        & CG-Descent 5.3   &   33     &   6.059e-015      &   4.697e-008            \\
\hline
deconvu      &  61    & mNewtow          &   84016   &   1.567e-009      &   9.999e-007            \\ 
             &        & L-CG-Descent     &   475     &   1.189e-008      &   9.187e-007            \\
             &        & L-BFGS           &   208     &   2.171e-010      &   8.924e-007            \\
             &        & CG-Descent 5.3   &   475     &   1.184e-008      &   9.078e-007            \\
\hline
denschna     &  2     & mNewtow          &   6     &   1.103e-023      &   6.642e-012            \\ 
             &        & L-CG-Descent     &   9     &   3.167e-016      &   3.527e-008            \\
             &        & L-BFGS           &   9     &   3.167e-016      &   3.527e-008            \\
             &        & CG-Descent 5.3   &   9     &   7.355e-016      &   4.825e-008            \\
\hline
denschnb     &  2     & mNewtow          &   6     &   5.550e-026      &   4.370e-013            \\ 
             &        & L-CG-Descent     &   7     &   3.641e-017      &   1.034e-008           \\
             &        & L-BFGS           &   7     &   3.641e-017      &   1.034e-008            \\
             &        & CG-Descent 5.3   &   8     &   4.702e-014      &   4.131e-007            \\
\hline
denschnc     &  2     & mNewtow          &   11     &  1.119e-021       &   1.731e-010            \\ 
             &        & L-CG-Descent     &   12     &  3.253e-019       &   3.276e-009            \\
             &        & L-BFGS           &   12     &  3.253e-019       &   3.276e-009            \\
             &        & CG-Descent 5.3   &   12     &  {\bf 1.834e-001} &   4.143e-007            \\
\hline
denschnd     &  3     & mNewtow          &   40     &  3.238e-010       &   9.897e-007            \\ 
             &        & L-CG-Descent     &   47     &  4.331e-010       &   8.483e-007            \\
             &        & L-BFGS           &   47     &  4.331e-010       &   8.483e-007            \\
             &        & CG-Descent 5.3   &   45     &  8.800e-009       &   6.115e-007            \\
\hline
denschnf     &  2     & mNewtow          &   6     &   6.513e-022      &    6.281e-010           \\ 
             &        & L-CG-Descent     &   8     &   2.126e-015      &    6.455e-007           \\
             &        & L-BFGS           &   8     &   2.126e-015      &    6.455e-007           \\
             &        & CG-Descent 5.3   &   11     &  1.104e-017      &    6.614e-008           \\
\hline
engval2      &  3     & mNewtow          &   13     &  2.199e-019       &   3.603e-008            \\ 
             &        & L-CG-Descent     &   26     &  1.034e-016       &   8.236e-007            \\
             &        & L-BFGS           &   26     &  1.034e-016       &   8.236e-007            \\
             &        & CG-Descent 5.3   &   76     &  3.185e-014       &   5.682e-007            \\
\hline
hairy        &  2     & mNewtow          &   19     &  2.000e+001       &   1.149e-008            \\ 
             &        & L-CG-Descent     &   36     &  2.000e+001       &   7.961e-011            \\
             &        & L-BFGS           &   36     &  2.000e+001       &   7.961e-011            \\
             &        & CG-Descent 5.3   &   14     &  2.000e+001       &   1.044e-007            \\
\hline
heart6ls     &  6     & mNewtow          &   312     &  1.038e-023       &   2.993e-008            \\ 
             &        & L-CG-Descent     &   684     &  2.646e-010       &   5.562e-007            \\
             &        & L-BFGS           &   684     &  2.646e-010       &   5.562e-007            \\
             &        & CG-Descent 5.3   &   2570    &  1.305e-010       &   2.421e-007            \\
\hline
helix        &  3     & mNewtow          &   13     &   3.585e-022      &    3.326e-010           \\ 
             &        & L-CG-Descent     &   23     &   1.604e-015      &    3.135e-007           \\
             &        & L-BFGS           &   23     &   1.604e-015      &    3.135e-007           \\
             &        & CG-Descent 5.3   &   44     &   2.427e-013      &    6.444e-007           \\
\hline
himmelbb     &  2     & mNewtow          &   7      &   7.783e-021      &    1.325e-007           \\ 
             &        & L-CG-Descent     &   10     &   9.294e-013      &    2.375e-007           \\
             &        & L-BFGS           &   10     &   9.294e-013      &    2.375e-007           \\
             &        & CG-Descent 5.3   &   11     &   1.584e-013      &    1.084e-008           \\
\hline
himmelbh     &  2     & mNewtow          &   4     &    -1.000e+000     &    1.085e-009           \\ 
             &        & L-CG-Descent     &   7     &    -1.000e+000     &    2.892e-011           \\
             &        & L-BFGS           &   7     &    -1.000e+000     &    2.892e-011           \\
             &        & CG-Descent 5.3   &   7     &    -1.000e+000     &    1.381e-007           \\
\hline
humps        &  2     & mNewtow          &   37     &    1.695e-013     &    1.826e-007           \\ 
             &        & L-CG-Descent     &   53     &    3.682e-012     &    8.552e-007           \\
             &        & L-BFGS           &   53     &    3.682e-012     &    8.552e-007           \\
             &        & CG-Descent 5.3   &   48     &    3.916e-012     &    8.774e-007           \\
\hline
jensmp       &  2     & mNewtow          &   10     &    1.244e+002     &    2.046e-012           \\ 
             &        & L-CG-Descent     &   15     &    1.244e+002     &    5.302e-010           \\
             &        & L-BFGS           &   15     &    1.244e+002     &    5.302e-010           \\
             &        & CG-Descent 5.3   &   13     &    1.244e+002     &    4.206e-009           \\
\hline
kowosb       &  4     & mNewtow          &   10     &    3.075e-004     &    1.055e-007           \\ 
             &        & L-CG-Descent     &   17     &    3.078e-004     &    3.704e-007           \\
             &        & L-BFGS           &   17     &    3.078e-004     &    3.704e-007           \\
             &        & CG-Descent 5.3   &   66     &    3.078e-004     &    8.818e-007           \\
\hline
loghairy     &  2     & mNewtow          &   23     &    1.823e-001     &    1.880e-007           \\ 
             &        & L-CG-Descent     &   27     &    1.823e-001     &    1.762e-007           \\
             &        & L-BFGS           &   27     &    1.823e-001     &    1.762e-007           \\
             &        & CG-Descent 5.3   &   46     &    1.823e-001     &    7.562e-008           \\
\hline
mancino      &  100   & mNewtow          &   5      &     1.257e-021    &     4.659e-008          \\ 
             &        & L-CG-Descent     &   11     &     9.245e-021    &     7.239e-008          \\
             &        & L-BFGS           &   9      &     3.048e-021    &     1.576e-007          \\
             &        & CG-Descent 5.3   &   11     &     9.245e-021    &     7.239e-008          \\
\hline
maratosb     &  2     & mNewtow          &   7        &   -1.000e+000      &   9.342e-011            \\ 
             &        & L-CG-Descent     &   1145     &   -1.000e+000      &   3.216e-007            \\
             &        & L-BFGS           &   1145     &   -1.000e+000      &   3.216e-007            \\
             &        & CG-Descent 5.3   &   946      &   -1.000e+000      &   3.230e-009            \\
\hline
mexhat       &  2     & mNewtow          &   4      &    -4.010e-002     &     1.972e-011          \\ 
             &        & L-CG-Descent     &   20     &    -4.001e-002     &     4.934e-009          \\
             &        & L-BFGS           &   20     &    -4.001e-002     &     4.934e-009          \\
             &        & CG-Descent 5.3   &   27     &    -4.001e-002     &     3.014e-007          \\
\hline
palmer1c     &  8     & mNewtow          &   7      &    9.760e-002     &     6.619e-007          \\ 
             &        & L-CG-Descent     &   11     &    9.761e-002     &     1.254e-009          \\
             &        & L-BFGS           &   11     &    9.761e-002     &     1.254e-009          \\
             &        & CG-Descent 5.3   &   126827 &    9.761e-002     &     9.545e-007          \\
\hline
palmer2c     &  8     & mNewtow          &   1      &   1.442e-002      &     1.023e-008          \\ 
             &        & L-CG-Descent     &   11     &   1.437e-002      &     1.257e-008          \\
             &        & L-BFGS           &   11     &   1.437e-002      &     1.257e-008          \\
             &        & CG-Descent 5.3   &   21362  &   1.437e-002      &     5.761e-007          \\
\hline
palmer3c     &  8     & mNewtow          &   1      &   1.954e-002      &     3.958e-009          \\ 
             &        & L-CG-Descent     &   11     &   1.954e-002      &     1.754e-010          \\
             &        & L-BFGS           &   11     &   1.954e-002      &     1.754e-010          \\
             &        & CG-Descent 5.3   &   5536   &   1.954e-002      &     9.753e-007          \\
\hline
palmer4c     &  8     & mNewtow          &   1      &   5.031e-002      &     1.123e-008          \\ 
             &        & L-CG-Descent     &   11     &   5.031e-002      &     3.928e-009          \\
             &        & L-BFGS           &   11     &   5.031e-002      &     3.928e-009          \\
             &        & CG-Descent 5.3   &   44211  &   5.031e-002      &     9.657e-007          \\
\hline
palmer5c     &  6     & mNewtow          &   1     &    2.128e+000     &      1.447e-013         \\ 
             &        & L-CG-Descent     &   6     &    2.128e+000     &      3.749e-012         \\
             &        & L-BFGS           &   6     &    2.128e+000     &      3.749e-012         \\
             &        & CG-Descent 5.3   &   6     &    2.128e+000     &      2.629e-009         \\
\hline
palmer6c     &  8     & mNewtow          &   1      &   1.639e-002      &     7.867e-010          \\ 
             &        & L-CG-Descent     &   11     &   1.639e-002      &     5.520e-009          \\
             &        & L-BFGS           &   11     &   1.639e-002      &     5.520e-009          \\
             &        & CG-Descent 5.3   &   14174  &   1.639e-002      &     7.738e-007          \\
\hline
palmer7c     &  8     & mNewtow          &   1      &   6.020e-001      &     9.090e-009          \\ 
             &        & L-CG-Descent     &   11     &   6.020e-001      &     7.132e-009          \\
             &        & L-BFGS           &   11     &   6.020e-001      &     7.132e-009          \\
             &        & CG-Descent 5.3   &   65294  &   6.020e-001      &     9.957e-007          \\
\hline
palmer8c     &  8     & mNewtow          &   1      &   1.598e-001      &     1.099e-009          \\ 
             &        & L-CG-Descent     &   11     &   1.598e-001      &     2.376e-009          \\
             &        & L-BFGS           &   11     &   1.598e-001      &     2.376e-009          \\
             &        & CG-Descent 5.3   &   8935   &   1.598e-001      &     9.394e-007          \\
\hline
rosenbr      &  2     & mNewtow          &   20     &   2.754e-013      &     8.253e-007          \\ 
             &        & L-CG-Descent     &   34     &   4.691e-018      &     7.167e-008          \\
             &        & L-BFGS           &   34     &   4.691e-018      &     7.167e-008          \\
             &        & CG-Descent 5.3   &   37     &   1.004e-014      &     1.894e-007          \\
\hline
sineval      &  2     & mNewtow          &   41     &   5.590e-033      &     2.069e-015          \\ 
             &        & L-CG-Descent     &   60     &   1.556e-023      &     1.817e-011          \\
             &        & L-BFGS           &   60     &   1.556e-023      &     1.817e-011          \\
             &        & CG-Descent 5.3   &   62     &   1.023e-012      &     5.575e-007          \\
\hline
sisser       &  2     & mNewtow          &   15    &    3.814e-010     &     4.485e-007          \\ 
             &        & L-CG-Descent     &   6     &    6.830e-012     &     2.220e-008          \\
             &        & L-BFGS           &   6     &    6.830e-012     &     2.220e-008          \\
             &        & CG-Descent 5.3   &   6     &    3.026e-014     &     3.663e-010          \\
\hline
tointqor     &  50    & mNewtow          &   1      &   1.176e+003      &     3.197e-014          \\ 
             &        & L-CG-Descent     &   29     &   1.175e+003      &     4.467e-007          \\
             &        & L-BFGS           &   28     &   1.175e+003      &     7.482e-007          \\
             &        & CG-Descent 5.3   &   29     &   1.175e+003      &     4.464e-007          \\
\hline
vardim       &  200   & mNewtow          &   19     &   1.365e-025      &     7.390e-011          \\ 
             &        & L-CG-Descent     &   10     &   4.168e-019      &     2.582e-007          \\
             &        & L-BFGS           &   7      &   5.890e-025      &     3.070e-010          \\
             &        & CG-Descent 5.3   &   10     &   4.168e-019      &     2.582e-007          \\
\hline
watson       &  12    & mNewtow          &   16     &   4.202e-006      &     1.918e-009          \\ 
             &        & L-CG-Descent     &   49     &   1.592e-007      &     8.026e-007          \\
             &        & L-BFGS           &   48     &   9.340e-008      &     1.319e-007          \\
             &        & CG-Descent 5.3   &   726    &   1.139e-007      &     8.115e-007          \\
\hline
yfitu        &  2     & mNewtow          &   37     &   6.670e-013      &     2.432e-012          \\ 
             &        & L-CG-Descent     &   75     &   8.074e-010      &     3.910e-007          \\
             &        & L-BFGS           &   75     &   8.074e-010      &     3.910e-007          \\
             &        & CG-Descent 5.3   &   147    &   2.969e-011      &     5.681e-007          \\
\hline
\end{longtable}
\end{center}
We summarize the comparison of the test results as follows:
\begin{itemize}
\item[1.] the modified Newton function {\tt mNewton} converges in all the test problems after
terminate condition $\| g(x_k) \|_{\infty} < 10^{-6}$ is met. For all problems except {\tt bard}, 
{\tt cliff}, {\tt deconvu}, {\tt sisser}, and {\tt vardim}, the modified Newton uses fewer iterations
than {\tt L-CG-Descent}, {\tt L-BFGS}, and {\tt CG-Descent 5.3} to converge to the minimum. Since
in each iteration, {\tt mNewton} needs more numerical operations, small iteration count does not mean
superior efficiency, but it indicates some promising.
\item[2.] For all the problems except the problem {\tt arglina}, all algorithms find the same mininum.
For the problem {\tt arglina}, the modified Newton finds a better local minimum.
\end{itemize}

Based on these test results, we believe that the new algorithm is promising. This leads us to use
the similar idea described in this paper to develop a modified BFGS algorithm which is very 
promising in the numerical tests.

\section{Conclusions}
We have proposed a modified Newton algorithm and proved that the modified Newton algorithm
is globally and quadratically convergent. We show that there is an efficient way to 
implement the proposed algorithm. We present some numerical test results. The results
show that the proposed algorithm is promising. The Matlab implementation {\tt mNewton} described in this 
paper is available from the author.

\section{acknoledgement}
The author would like to thank Dr. Sven Leyffer at Atgonne National Laboratory
for his suggestion of comparing the performance of 
Algorithm \ref{implementedAlgo} to some established and state-of-the-art algorithms.
The author is also grateful to Professor Teresa Monterio at Universidade do Minho for
providing the software to convert AMPL mod-files into nl-files, which makes the test possible.

%\bibliography{Myrefs}  % The name of your .bib file

\bibliographystyle{siam}     % The style of your bibliography

\end{document}